\documentclass[12pt]{amsart}
\usepackage{amsmath,amssymb,amsfonts,amscd, graphicx, latexsym, verbatim, multirow, color}
\usepackage{amsthm}
\usepackage{tikz}
\usepackage{forest}
\usepackage{float}
\usepackage{listings}
\usepackage{adjustbox}
\usepackage{arabtex}
\usetikzlibrary{decorations.shapes}
\theoremstyle{definition}

\makeatletter
\@addtoreset{equation}{section}
\makeatother

\providecommand{\AMS}{$\mathcal{A}$\kern-.1667em%
\lower.25em\hbox{$\mathcal{M}$}\kern-.125em$\mathcal{S}$}

\usepackage{courier}

\begin{document}
\title[Caf\'e conversations]{Caf\'e conversations: \\ A tribute to Norbert A'Campo}

\author{Athanase Papadopoulos}
\address{Athanase Papadopoulos,  Universit{\'e} de Strasbourg and CNRS,
7 rue Ren\'e Descartes,
 67084 Strasbourg Cedex, France}
\email{papadop@math.unistra.fr}
 
\date{}
\maketitle

\noindent \emph{Abstract.} This is an intrusion in the life and the mathematics of Norbert A'Campo, intended to be a tribute to him and an acknowledgement of his impact on those who know him and his work.
\\
The final version of this paper appears in the book ``Essays in Geometry, Dedicated to Norbert A’Campo",  IRMA Lectures in Mathematics and Theoretical Physics, EMS Press, Berlin, 2023.
\bigskip

AMS codes: 01-02, 01-06, 01A70

Keywords: Biography, Norbert A'Campo, mathematics in France.

\bigskip

 \bigskip

 -- ``Excuse me, what are you talking about ?"
 
 \medskip

 We have heard this question so many times, from onlookers who listened to us talk for hours in caf\'es in Basel (Caf\'e Euler), Strasbourg (Christian and Caf\'e Brant), Vienna (Caf\'e Sperl), Paris (Rostand), other small or grand caf\'es, sometimes railway station caf\'es, in Mulhouse, Djursholm, Trento, Istanbul, Cagliari, Karlovassi, Mandraki, Barcelona, Beijing, Bangalore, Delhi, Varanasi, and other cities scattered in our Eurasian continent where mathematics drove us.  
Sometimes, the waiters themselves, if they don't know us yet, ask this same question; they are curious to know what is this gibberish that we have been uttering for hours. Occasionally, we engage in a conversation with them. Usually, Norbert gives the answer: ``We are drawing dessins d'enfants", or: ``We are discussing infinity", or: ``We are counting singular points", or: ``We are considering the field with one element", or: ``We are constructing branched coverings", or: ``We are counting pairs of pants", and other phrases like these which put even more mystery in the eyes of the stunned interlocutor. To a waiter at the caf\'e of the Grand Hotel Euler in Basel, Norbert once asked: ``Do you know who is the Grand Euler?" The waiter had no idea. That is how we occasionally made friends at the caf\'e. 
 
``Mathematics is made for travelling" says Norbert, a sentence we often repeat to each other. Physical journeys, and also journeys in the world of ideas. When I discuss with him, it is always a journey.

It was during those countless conversations at the caf\'e that I learned much of what I know in mathematics: spherical geometry, the Riemann--Roch theorem, Serre duality, Dolbeault's lemma, transitional geometries, the importance of the one-element field, Belyi's theorem, slalom curves and polynomials, and many other things. 
I also learned some details about Norbert's life story and that is what I will talk about in the next few pages.

\medskip

Norbert was born just over 80 years ago in a farm that belonged to his family. His father together with his uncle took it over from their parents. When Norbert was born, the extended family lived there and everybody took part in the manifold farm duties. 

The farm was situated in Merkelbeek, a small village in the province of Limburg, five kilometers from the German border. The province has been part of the Netherlands since 1815.  
Today, the former village structure has disappeared due to a reorganization of the region by the central government.

A'Campo is a Latin-sounding name. It is conceivable that the origin is Van der Velde. In French it would have been Deschamps, in English Fields, and in German vom Feld. 
The name was probably latinized around 1850, because of the ruling power: people did not want their name to sound too Dutch, nor too German.

Norbert's mother's ancestors were Germans who had emigrated to the Netherlands to escape Bismarck. She was also raised in a farm, and after  her marriage she identified herself completely with the new farm life.

Back when Norbert was a child, in order to buy items that were a little different from the everyday things (shoes, farm equipment, etc.), the A'Campo family used to go to Heerlen, a name which means ``brighter": you can see the sun earlier there. Indeed, Heerlen is situated in the east of the Limburg. It is a small town with a hospital and a veterinarian, merchants of agricultural equipment and people who could repair tractor wheels.

 The nearest big city was Maastricht. Among the other cities that people from Merkelbeek could visit easily were Aachen, Cologne and Li\`ege.

Norbert couldn't stand the Kindergarten. After a few days of trial and refusal, he was allowed to stay at home. The farm was quite big.   One has to imagine that there were many children there and wonderful opportunities for them to play together, much more than at the Kindergarten. At the same time, Norbert could observe his father doing his numerous agricultural tasks. It is in this environment where man has a deep relationship with the Earth that his mind was slowly shaped. ``\emph{O fortunatos nimium, sua si bona norint, agricolas!}" wrote Virgil about the charm of the country life.
 

For elementary school, Norbert went to Brunssum, the nearby village. The studies there lasted six years. At the beginning, a young girl came every day to pick him up on her bicycle and take him to school. Later on, he started to go alone, on foot. From age 8 or 9, he was permitted---although he did this rarely---to go to school on horseback.  Among the several horses in the farm, the one he was allowed to ride was used for light duties such as small transport: riding or pulling an open carriage. At times, the family used such a carriage to visit the grandmother who lived in another village, about 20 km away.  The other horses
were used for hard work (ploughing, etc.). Starting in the fifties, working horses were gradually replaced by tractors.

At school, when Norbert arrived on  horseback, he used to leave the horse in a pasture that belonged to a friend of his father, located in front of the school.  
He liked to keep an eye on it while he was in class, and  for that purpose he would sit by the window. Although it was rare that he went at school on horseback, his school records of those times read: ``Does not always pay attention in class," or  ``Often looks out the window".

Norbert remembers that in the mid 1950s, the Russians, in their will to develop agriculture, came to the Netherlands to buy cows. His father sold a number of them. To take care of them during the weeks-long travel on the train was not easy, and Norbert asked his father to allow him to accompany the cows aboard the same train to help with this task. His mother did not agree.

At the school in Brunssum, Norbert used to take additional lessons with a private teacher---his father wanted him to learn High Dutch---, but actually the lessons ended up being math lessons because the teacher realized that Norbert liked it. He introduced him very early to calculus with the variables $x$ and $ij$  (in the  Dutch system, the second variable is usually not called $y$, but $ij$).

After the six years of elementary school, Norbert was put as a boarder at the secondary school in Weert, 60 km from Merkelbeek.  This distance was quite large in those days. 
The duration of the studies at the secondary school was also six years. It was an excellent classical high school with several teachers who had completed advanced studies. The Dutch teacher spoke several languages. He had a theory about language development and he was particularly interested in Norbert's native Limburg language. At the same time, there was a very solid language program in that school: Latin, Greek, German, French, English, and of course Dutch. There was very little mathematics and science. 
 Norbert learned at school how to debate on something. Indeed, in the Dutch school, pupils learned the art of rhetoric, that is, how to defend a thesis, even if you do not agree with it. The pupils were taught how to maintain the upper position in a discussion.
  This was part of the youth education.

The young Norbert thought he would become a farmer later on. He used to watch his father working, and he liked that work. Had he made that choice, a farm was ready for him. The entire surrounding, that is, the extended family together with the employed people, was centered on the farm work. In that society, it was known that eventually the young become stronger than the older ones, and then the farm starts to belong to them. Some farmers suffer from this, and sometimes fall into depression. Today, Norbert's brother has taken over the farm. All his life, Norbert felt certain that he could always return to the farm as an alternative to academic life.

Many things have changed in the region since that time, including the language, which became a  rarely spoken dialect. But very recently there has been a will to recover the old language there. In the past, it was said that this language sounded too much like German. Today, the local government is trying to introduce it again  in primary education. Norbert told me several times that there are similarities between this dialect and the Basel  dialect which is spoken in the village he lives in presently. He has a theory that says that there are common features to all the languages spoken by people living along the Rhine.

From his years at the farm, Norbert has kept the sense of country life. He is the opposite of a city dweller.  
Later, in the Paris region, when he got a job at Jussieu (University of Paris VII), then Orsay (University of Paris XI), he lived in Chantecoq, a name that fits well with the place.  In Basel, the A'Campo family lives in Witterswil, a small village in the countryside. 

It was somewhat by the luck of an encounter that Norbert decided to study mathematics. 
The story  starts with the first known relative who actually went to university, an uncle of Norbert's father, called Emile.   The latter, in his youth, wanted to study mathematics, but his parents did not accept. They agreed however that he would study biology, which seemed to them to be a relatively useful field. Thus, Uncle Emile entered university and studied biology. After that, he went into business and became a very successful entrepreneur. At age 45, he began taking advanced mathematics courses, just for pleasure. He attended Emile Borel's lessons in Paris. When he went to Paris, he used to stay at the George V. This was unusual for a mathematics student!
 
Uncle Emile visited the A'Campo family from time to time, and he used to go there with his car.  At that time, very few people they knew had a car;  in the whole village there were only two cars. 
Norbert remembers that in elementary school, children used to enjoy counting the number of vehicles that passed by. At the end of the day, there were five or so.
 Norbert's family didn't have a car; they used the horse.  A car, in those times, was linked to a function: people rented one, with driver, only when they had to go somewhere.  If the trip was to last three days, the driver would charge a reduced price. Norbert's father bought his first
car around 1956-57.

 At secondary school, Norbert liked physics and initially he thought that later he would study physics at university.  
Uncle Emile always spoke of the profession of mathematician as an extraordinary profession. 
 He had told him once: ``If you want to enter university, I have a friend who is a mathematician who can advise you." This friend was a well-known professor of mathematics, Jaap Seidel, teaching in Eindhoven.  Norbert and his uncle went to see him. Seidel said to Norbert: ``You want to do physics, that's fine, but you should know that in this field there is too much of a hierarchy. To do research in physics, you need a laboratory and a lot of money. If you become a physicist, it is only at the age of 50 that you can become autonomous. In mathematics it is not the same; you are free and independent since the beginning."  He added: ``If you decide to become a mathematician, the best place to study math is Utrecht."
 
 Uncle Emile also told Norbert  about his mathematical dream, namely to find an explanation of how living cells or more generally life is possible contrary to the growing entropy, by finding an extra term to insert in the Boltzmann equation which (for a limited time) can slow down or even reverse the growth of entropy.

It was under the influence of this uncle that Norbert turned to mathematics.
His father  was always very supportive and he encouraged him in that direction.

Thus, Norbert enrolled in mathematics at the University of Utrecht. 

 Among the professors who marked him there, I heard him mentioning several times the names Freudenthal, 
van der Blij and de Iongh.

  Norbert remembers Freudenthal's first lecture: ``Consider a plane, such as the blackboard plane, divided into two parts by a vertical line, where the straight lines change their direction when crossing this separating line, such as rays of light entering water." Freudenthal's suggestion was to study the geometry of this kind of plane with a new concept of line, talking about axioms like the existence of a unique line through any given pair of points, etc.
  
Norbert also remembers the proof that Freudenthal gave them of the fact that a set cannot be in one-to-one correspondence with the set of its subsets. I didn't ask him if he remembered other details, but I can imagine
 that Freudenthal took this opportunity to tell them about the continuum hypothesis and the important developments it led to.

 Hans Freudenthal had moved to the Netherlands after fleeing Nazi Germany.  He had been an assistant to Brouwer in Amsterdam.  As early as the first year, he explained to his students how to do math. He advised them to try to understand the lectures without taking notes. At the end of his lectures, he used to ask the students if they found it complicated. He would also tell them: ``Next time, one of you will try to explain to the whole class what I taught you today." For him, that was the way to do progress in math:  to explain without looking at notes. Usually no student could do it right the first time. They  would start over.
 
Norbert kept these principles alive. I have seen him attending conferences, and giving lectures, hundreds of times, over several decades. I never saw him taking notes or looking at notes.

The amount of knowledge that Freudenthal gave to his students was minimal but profound. A specialist in topology, he was no less a specialist in history of mathematics and more generally in history of science. In his teaching, he always tried to relate mathematical problems to those of everyday life, a point of view that Norbert always adopted later on.

Norbert also remembers a very good analysis course by Frederik van der Blij. The latter, like Freudenthal, gave him a taste for  history of mathematics and the connections between mathematics and everyday life.
With van der Blij, the atmosphere was comparable to that with Freudenthal. 
He used to explain to them profound notions by doing simple things. For instance: Consider the sine formula given by its power series expansion. Try to find where the zeros are. Try to prove that there is a smallest positive zero; that would be $\pi$.  

From de J. J. de Iongh, Norbert remembers his seminars on set theory and  G\"odel's constructions.

At the University of Utrecht, people sometimes talked about French mathematics. Everybody knew that if you wanted to become a  mathematician, it was advisable to go to France, a country where mathematics was highly valued. 
Some weekends, the mathematics assistants at the university of Utrecht used to come back very tired from Paris, but they all used to say: it was fantastic; they were coming back from the Bourbaki seminar.

Norbert's father had a friend, an old classmate, who was living in N\^\i mes. This friend had joined the resistance in France during the Second World War. After the war he remained in France and became a construction equipment dealer in N\^\i mes. 
One of Norbert's father's dreams was to see this friend again. Eventually, he went to visit him.  At his return from France, he gave a very picturesque description of the countryside in N\^\i mes.  Norbert, who at that time was a student in Utrecht, was fascinated. He inquired about N\^\i mes. He found out that there was no university there, but that there was one in Montpellier, 60 km away from N\^\i mes.  He wrote to the university of Montpellier, asking whether he could continue his mathematics education there. The answer came quickly, a very gentle letter saying that this was possible.  In those days, it was not complicated to be accepted at a French university. 
This is how Norbert became a student at the University of Montpellier. 

Norbert arrived in Montpellier towards the end of 1961. The 1961-62 academic year had already begun, and he enrolled there in the mid of that  academic year, in January 1962. Since he had already completed two years in Utrecht, he was admitted directky to the third year.

In Montpellier, several circumstances played in his favor.
The first one took place during the second year of his studies there. The teaching assistant of the course on holomorphic functions fell ill. In the meanwhile, Norbert had done all the exercises. The course coordinator, who knew that, suggested that Norbert replace the assistant.  The job consisted in helping the students with the exercises. By doing this, Norbert went more profoundly in the subject of holomorphic functions. At the same time, this allowed him to earn some money.

  Another event marked Norbert during the same period. 
 At the end of the autumn of 1962, about a year after his arrival in Montpellier, he went on a grape harvesting trip to Vendargues, a village situated about 10 km from Montpellier. For him, it was also a way of returning to the soil, of helping to collect from it the gifts it offers to man.
The winegrower for which he was working asked him what he was doing in life. Norbert  responded that he was studying mathematics. The winegrower then asked him if he knew Grothendieck.  Norbert said  no, I don't know him; he had never heard that name.  

It was only about two years later that Norbert heard the name again, when Jean-Pierre Lafon, one of his teachers in Montpellier, asked him to read Grothendieck's Tohoku paper,  \emph{On some
points of homological algebra}.  Norbert then remembered the winegrower who had pronounced this Dutch-sounding name. He was puzzled and he returned to Vendargues, where he 
 asked the peasant how come he knows Grothendieck. The winegrower answered: ``Everyone here knows Grothendieck. He spent his childhood with his mother  in a house in Meyrargues, at a short distance from
Vendargues." When the winegrower saw that Norbert was very much interested, he took him by tractor to show him the house. 

  Grothendieck had completed high school and university in Montpellier.  At the time we are talking about, he was making his career in Paris. Later, he returned again to Montpellier.

Norbert still remembers the name of the winegrower who first told him about Grothendieck, Jean-Henri Teissier. He also remembers his wife Annette, their son Jean-Louis and their daughter, Yvette.   He returned later to Vendargues. He wanted to see again Grothendieck's house in Meyrargues, but he could not find his way. When he went there for the first time with Teissier, they were chatting together on the tractor, and he did not notice the way. In those times there were no paved roads and the tractor went through the fields. 
 Today there is a national road there, and even, the highway is close by. It is still difficult for Norbert to find his way around, because Vendargues has completely changed. 
 
 Norbert returned several times to visit Teissier, with whom he became friends. Teissier no longer lives. Later, with his wife Annette (the same name), Norbert visited the farm again. They talked with Teissier's wife and their son. The latter studied geology in Bordeaux. Looking in my old mails, I find the following, from Norbert, dated May 23 mai 2016:
``Dear Athanase, 
  I came back from Montpellier. 
 I went to Vendargues and Meyrargues.
Grothendieck's footsteps are now very hidden.  Meyrargues did not change. 
But Vendargues became a little city. I bought several bottles of wine from Meyrargues."

In the 1970s, after Grothendieck withdrew from the Parisian world  of mathematics and returned living in a commune near Montpellier, he started with some students at the University of Montpellier   a program on Thurston's theory of surfaces, in which he involved Norbert. 

Grothendieck and Norbert shared the same worries about the preservation of nature and of the animal species---including the humans. 
Grothendieck had decided to make from the return to nature a political commitment, but that is another story.  My friend and collaborator Stelios Negrepontis, who was teaching at McGill at that time, went to Nice for the 1970 ICM. He remembers seeing Grothendieck talking there about \emph{Survivre}, the antimilitarist and ecologist movement he had founded,  together with  Claude Chevalley and Pierre Samuel, a couple of years before. Stelios, who had an awe for Grothendieck, told me that he was surprised by the somewhat insolent manner that
young French mathematicians were addressing him, not at all convinced by the value of his Survivre action. 

Let us return to Thurston's theory. 

Grothendieck's program included the classification, up to isotopy and mapping class group action, of objects which could be systems of curves, subsurfaces and mixtures of them. A group appeared, the cartographic group, which is related to dessins d'enfants and to other classes of cell decompositions of a surface. There is also a universal cartographic group. In working on this topic, Grothendieck developed a sort of combinatorial Teichm\"uller theory. One motivation for his interest in this subject was the Nielsen realization problem. This was before Steve Kerckhoff solved the problem. The solution came amid this Montpellier activity, and Grothendieck asked Norbert to explain to him and his students Kerckhoff's work. At the same time, he asked him to be a member of the jury committee of the doctoral dissertation of his student Yves Ladegaillerie. Later, Norbert became the formal advisor of another  student  of  Grothendieck, Pierre Damphousse. The latter's thesis defense took place at Orsay, in June 1981. Norbert is not very much talented for administrative work, but Grothendieck, at that time,  was even less. 
  I met Damphousse  at the 2010 Clay conference celebrating the proof of the Poincar\'e conjecture. It was at Thurston's talk which was held in the magnificent lecture hall of Oceanographic Institute. I was there with Norbert, who introduced me to Damphousse. The latter died unexpectedly shortly after, at a relatively young age.

At several occasions,  I asked Norbert who were the mathematicians who influenced him. I heard the name Marcel Lefranc, an excellent mathematician, he said. Norbert followed his lectures on homology.  Thanks to these lectures,   he learned the  Lefschetz fixed point  formula.  He also learned there Sperner's lemma, which gives an alternative proof of Brouwer's fixed point theorem, without using the notion of degree. In the same period, Norbert purchased Spanier's book to learn more topology. But it was through Lefranc's course that he understood the idea of homology. Norbert also told me he learned a lot from Roger Cuppens, who was the teaching assistant of the probability course. He remembers that the latter explained to him the content of his thesis on infinitely divisible laws. Later, Cuppens became a professor in Toulouse.

At the time where Norbert was studying in Montpellier, there was a very dense mathematical activity tthere, whose most important part---according to Norbert---took place at a caf\'e, ``Chez Jules", in front of the faculty of mathematics, which was still located on rue de l'Universit\'e, in front of the main university building.  Today, the Faculty of Sciences of the University of Montpellier, including the mathematics department,  is situated out of the historical center; only the administration   remains in the original building.

Among the other professors in Montpellier to whom Norbert owes part of his mathematical education, he mentioned to me several times the name of Andr\'e Martineau, who taught holomorphic functions in several variables. Martineau liked to explain to his students works of other mathematicians, and 
Norbert learned from him
a lot of mathematics. Sometimes Martineau invited the students to his home and continued to talk there while he looked after his children. His wife was a literature teacher at the University of Montpellier. Norbert learned from Martineau   the Levy problem, holomorphic domains and the Cartan--Oka theorem.  Martineau did a lot of mathematics with the students at the Caf\'e Jules.

 Norbert also mentioned Lafon, who has obtained an analogue of the Weierstrass preparation theorem for certain subrings of the ring of convergent power series. Lafon asked then whether such a theorem remains valid for a certain kind of larger ring containing convergent power series, namely, a ring of  power series satisfying inequalities on the coefficients such as Gevrey classes. One day Norbert met Lafon at the beach and the latter told him about this question. Soon after, Norbert came back with a proof of the fact that the answer to this question is negative. Lafon told him then: ``You should write this proof and send it to Malgrange''. This is how Norbert got in touch with Malgrange. He wrote down his result and he asked Martineau whether this could give him a doctorate. Martineau told him: ``Leave it now in a drawer, and wait for something better."   A posteriori, Norbert thinks it was a good advice, not to make a PhD too early.

  At the beginning of the Academic year 1967--68, Norbert was hired at the University of Poitiers, as an assistant, working in the group of Pierre Dolbeault. By then, he had a car and he used it for moving  from Montpellier to Poitiers. In the meanwhile, he leaned about the existence of a conference on foliations at 
Mont Aigoual, a place on the way between Montpellier (in the South East of France) and Poitiers (in the South West). Mont Aigoual is an elevation, situated in the Massif Central, a region in the Central South of France, consisting of mountains and plateaus. It is also a resort place, with snow and skying.

Norbert likes snow. I remember a scene. I was coming back with him, in his car, from Italy to France and Switzerland.  It was in March 2016, we were heading towards Strasbourg where he proposed to drop me, before going to his home in Basel.  We had spent two weeks at the Mathematical Center  in Trento,  with Sumio Yamada, working on a project which we did not finish until today.
On our way back, we took the pass called the Julier Pass, which had just opened, that year, after the period of big snow. When we arrived up there,  Norbert stopped the car, went out and walked on the snow, just a few minutes. He touched the snow with his hands. I stayed in the car, watching him.

Back to 1967-68. 

During his trip from Montpellier to Poitiers, Norbert stopped in mont Aigoual. 
He heard there a talk by Harold Rosenberg, in which the latter mentioned the problem of foliating the 5-sphere. For Norbert, foliations was a new subject. 
Rosenberg explained the construction of Reeb's foliation of the 3-sphere and he said that the existence of a foliation on the 5-sphere was an open question.

Soon after, Norbert constructed a foliation of the 5-sphere. He went to Rosenberg, who was teaching at Orsay, and he showed him his foliation. Rosenberg said that this was good. Norbert asked him whether that could constitute a  PhD thesis, and Rosenberg said yes. 
At about the same time, Norbert proved that any simply connected 5-dimensional manifold whose second Stiefel--Whitney class is zero carries a codimension-1 foliation.

There was a lot of paper work to do for the thesis, and this took some time. In the meanwhile Norbert  proved also his monodromy theorem, 
which also became part of the thesis. Partly thanks to Milnor's book on isolated hypersurface singularities, monodromy had become an interesting topic for topologists. An important result of Grothendieck states that the homological eigenvalues of the monodromy are all roots of unity and the length of the Jordan blocks is bounded by the dimension. Brieskorn conjectured that the monodromy has no Jordan blocks and hence is of finite order.
Norbert constructed the first example of an isolated singularity with infinite homological local monodromy. This gave a counterexample to Brieskorn's conjecture and it became almost immediately an ingredient for Deligne's proof of the Weil conjecture. Directly afterwards, Malgrange produced similar examples in higher dimensions with Jordan blocks of the maximal possible length.  Today, Norbert is still working on monodromy.

 In Poitiers, like in Montpellier, Norbert learned a lot of mathematics. At that time Andr\'e Revuz was the head of the mathematics department there. Remarkably, he encouraged all the young members of the department working on their theses to travel to any conference or seminar they liked. He used to tell them that the department will pay their expenses provided that, when they come back, they give a seminar talk on what they heard. Those talks were usually very interesting, especially because Revuz was always asking for examples.

Pierre Dolbeault, who was also in Poitiers,  was a specialist of complex analysis of several variables. He had obtained his doctorate with Henri Cartan. 
He used to organize a seminar. One day, he was looking for someone who could give a talk.  Joseph Le Potier, who had just arrived from Rennes and who was still unknown, proposed to give a talk. He ended up talking for the whole semester. He explained the  Atiyah--Singer theorem with many of its applications. This was in 1968. 

This was the atmosphere in Poitiers at the time when Norbert was trained there.

A few years later, Le Potier, with Dolbeault as advisor, presented a Doctorat d'\'Etat in which he proved a conjecture of Griffith. After Poitiers, he obtained, like Revuz, a position at Orsay and then at Jussieu.
 
In mont Aigoual, Norbert had another favorable experience. Bill Browder was there and did not go to any of the talks. He was spending his time walking and visiting around, and he was happy to find someone like Norbert who would do the same.  They used to walk together, with Browder telling math to Norbert.  This is how Norbert learned about the $h$-cobordism theorem. He came back to Poitiers and recounted this to Dolbeault, who said that this was interesting, and he asked Norbert to give a course on this subject.
Norbert started a one-semester course on a theory  that he had heard for the first time a week before.  This is how he learned $h$-cobordism theory. Among the students at that course were Le Potier, Cathelineau, Poly, and many others.

From Poitiers, Norbert attended seminars in other cities.  He went regularly to Paris. In the train, he used to meet people going from Bordeaux to Paris like him to attend mathematical seminars. He liked to go the Institut Henri Poincar\'e (IHP), to listen to the talks of the number theory seminar organized by Delange, Poitou and Pisot. In 1967/68 and 1968/69 he  gave several talks there, on analogues of Malgrange's differentiable preparation theorem in the setting  of ultrametric valuated fields of arbitrary characteristic.
Krasner had also a seminar at IHP.
  The caf\'e \emph{La Contrescarpe}, a few hundred meters from IHP, was a meeting place, with Rosenberg and others going there regularly.

In the beginning of the seventies, the discussions following the regular seminars and colloquia in Paris were very lively and were often continued afterwards at some mathematicians' homes, several of which situated in upscale areas of Paris. Alain Chenciner lived in a studio apartment on \^Ile Saint Louis, arguably the prettiest part of the center of Paris. Sandy Blank lived right next door, Georges Pompidou as well. Together with other mathematicians such as Laudenbach, Norbert would spend hours there discussing mathematics.

After his first meeting with Browder, and after he became familiar with $h$-cobordism theory, Norbert wanted to learn more topology. Malgrange told him that for that he  should go to Orsay.  The University of Paris-Sud at Orsay was newly founded, with an ``\'Equipe de topologie" (topology team) working  around Cerf, who was Henri Cartan's student.  Some North-Americans, including Siebenmann and Rosenberg, were part of the team. Melvin Rothenberg was an invited professor. Norbert told me once that  the latter  was a very good lecturer. The young mathematicians in the Orsay team included Vogel, Latour, Barge, Lannes, Laudenbach and some others. Young researchers from Dijon working on foliations  (Joubert, Moussu, etc.) also participated in the Orsay topology seminar.   Ibish used to come quite regularly  from Nantes. Eventually, in Dijon and Nantes, they hired young people who were trained at Orsay, and gradually  these two cities became also centers for topology. 

 Lafon, who had moved to Toulouse, told Norbert that there was a conference there.  The reader might remember that in those times there were not as many conferences as there are today. Norbert went to that conference. Soon after the lectures started, someone said that there were interesting things happening in the next door building. They all went there. There were students, explaining their claims and demands.  We were in May 68.

In 1968-69, Norbert followed Thom's seminar on dynamics in Bures-sur-Yvette. Smale, Shub, Peixoto and others participated to that seminar. Norbert had already met Thom in Montpellier. He had attended there a talk by him in 1967.  He remembers that the latter had asked a question about the typical patterns caused by frost on windows.

 Norbert told me that at the same period, a change took place in the agricultural world. Before 1968, there were no big tractors in French farms,  and a substantial part of the work was done by hand. Famers needed a certain number of  workers, and these workers were poorly paid. The year 1968 saw the appearance of big agricultural machines, and Norbert realized this transformation. Each year, he used to go to the \emph{Salon international de l'Agriculture et de la machine agricole}, which was held in Paris, Porte de Versailles. He wanted  to see the new agricultural machines. His  father was particularly interested in how machines work, and Norbert learnt this curiosity from him. 
 In 1969, he saw there for the first time big tractors. Today, in the French farms, there are only big tractors, and very few people working.

 Back to Poitiers.

Pierre Dolbeault, like Norbert, came from a family of farmers.   In 1969, they went together to a conference on foliations in Oberwolfach. This was after Norbert tried to explain to Dolbeault his construction of a foliation of the five-dimensional sphere. Dolbeault was not familiar with the subject and he told Norbert: There is a conference on foliations in Oberwolfach, we can go together. 
 They used Norbert's car. During this trip, Norbert saw Dolbeault's parents' farm, which is located in central France, halfway between Poitiers and Oberwolfach. They spent a night there. The car needed some repair and Norbert repaired it there. In the farm, he could find everything that he needed for the repair. 

It was the first time that Norbert visited Oberwolfach.   The lectures were still held in the old building. At that conference, he met for the first time Reinhold Baer, Blaine Lawson, Robert Lutz, Robert Roussarie, Claude Godbillon, Jean Martinet, Georges Reeb, Jacques Vey, Andr\'e Haefliger and other mathematicians. Norbert explained his construction of a foliation of $S^5$, and Haefliger was very much interested. He asked him to explain to him all the details. Nobert tells me that he read very carefully what he had written, discussing every word, every comma. Between Norbert and Haefliger, it was the beginning of a lifelong friendship. It was also the starting point of Norbert's relation with Switzerland.

Norbert obtained his doctorate at the University of Paris-Sud, in Orsay, in 1972. The title of the dissertation was \emph{Feuilletages de vari\'et\'es et monodromie des singularit\'es} (Foliations of manifolds and monodromy of singularities).   At that time, the Orsay department was very young. It had been founded in  1965, first as part of the University of Paris, integrated in a project to decentralize this university. In 1971 it had  become a department of the Facult\'e des Sciences d'Orsay, in the newly founded Universit\'e de Paris-Sud. The creation of this university was an important step of the ungoing vast reorganization of Parisian universities, whose goal was to meet the demand for a transformation of the higher French educational system. Several mathematicians, including Henri Cartan, Hubert Delange and Georges Poitou moved from Paris to Orsay. 
Claude N\'eron, Jean Cerf, Adrien Douady, Michal Raynaud, Michel Demazure, Larry Siebenmann, Yves Meyer and Valentin Po\'enaru were also hired at Orsay.

   1972 is also the year when Thurston obtained his PhD, whose subject was also foliations. Thurston's dissertation was the starting point of a series of papers in which he solved all the important existing problems in the subject. The same year, Thurston visited Orsay. Norbert remembers the circumstances of that visit.  An invitation  was made to Morris Hirsch, by Rosenberg and Siebenmann, to visit Orsay. The paperwork was done, but eventually Hirsch declined, and the reason he gave was that he had a very talented student of which he had to take care. This student was Thurston. Rosenberg and Siebenmann decided then to invite the student; this is how Thurston arrived to Orsay.
  He lectured there on his version of what became known later as  the $h$-principle for foliations of codimension greater than 1. He also explained his result saying that an arbitrary field of 2-planes on a manifold of dimension at least four is homotopic to a field tangent to a foliation. Several young mathematicians working on foliations attended his lectures. During the same visit, Thurston went to Dijon and Plans-sur-Bex, and a large portion of these trips were done in Norbert's car. The latter remembers that at that time, Thurston was already thinking about hyperbolic geometry in dimension two. He asked  him how he came to know about this subject, and Thurston answered that it was his father who first told him about it.

   At that time, Norbert had a position of ``Ma\^\i tre de conf\'erences provisoire" at Orsay. Before that, he had taught for one year at Jussieu, also as  ``Ma\^\i tre de conf\'erences provisoire," for one year, 1970-71.
    
    Norbert's first visit to the United States took place in 1971, when Sandy Blank invited him to Boston. On that occasion, Blank arranged Norbert's first meeting with Dennis Sullivan. He  had told Norbert in advance that on the day of his arrival he will introduce him to Sullivan. Norbert had read some of the latter's works while he was a student in Poitiers.  In fact, 
  part of his construction of a foliation of the 5-sphere was inspired by Sullivan's work; in particular he had read his preprint on synthetic homotopy, which  contains a special decomposition of the 5-sphere. This decomposition is not as good as that of the 3-sphere which is used in the construction of the Reeb foliation, but with a little bit of thinking it  became the starting point for the construction of Norbert's foliation of  $S^5$. Only a slight modification was needed. 
  
 Sandy Blank organized an outing with Norbert and Dennis to attend a soccer or American football game (Norbert does not remember which of the two).  It was the last round of the American competition. The three of them went to the stadium. When they arrived, Norbert and Dennis had already started talking math. They immediately moved to the basement of the stadium, where they continued their discussion. Suddenly they realized that everybody in the stadium was leaving. The game was over, and they hadn't seen anything of it. They had forgotten about the game and why they came there. This is why Norbert does not know  whether it was soccer or American football: he did not see anything of the game. 

During his 1971-visit to the US, Norbert attended a course that Sullivan gave in Boston.
When he returned to Orsay, he told people about Sullivan and  his course.  Shortly after, there was an opportunity of inviting a mathematician as a one-year visiting professor at Orsay, and they decided it would be Sullivan. They needed to fill in some administrative documents, and they had only one day for that.  Cerf came to Norbert and asked for help in filling in the documents. Like others in Orsay, he considered that Norbert was close enough to Sullivan to know what to write on that form.  They had to specify where Sullivan was born and similar things. One has to remember that in those days email did not exist, and that it was not easy to communicate rapidly with someone in the United States. Cerf thought it is important to fill in the form correctly whereas Norbert was not so concerned about such details. To simplify the matter, he told Cerf that Dennis was born in Dallas. In fact, Norbert just thought that he was born in Texas, and the only city in Texas he could remember was Dallas. Thus, on Sullivan's official documents at Orsay it is written that he was born in Dallas. In fact, Sullivan was not born in Texas, but in Michigan.

I asked Sullivan what he remembers about this hiring process at Orsay. He twrote: ``Norbert had to do quickly. He  made up a fictitious address in Dallas, Texas, where JFK was assassinated ten years before. The first day at work at Orsay, Norbert introduced me to Michel Herman, Fran\c cois Laudenbach and Harold Rosenberg (at the Orsay pool in the latter case).  His family and mine spent several `promenade du dimanche apr\`es le repas'  days together. His daughter and mine later became friends as adults."

 Norbert did not like the division of the Orsay mathematics department into  ``\'Equipes". He felt that this compartmentalized people, and that the seminars had become more specialized. He, personally, did not consider himself more of a topologist than a geometer or algebraist. He remembers Douady saying: ``I don't mind dividing the department into teams, but then I want to be part of every team".

%
%

Norbert gave his first course at Orsay in 1973. 
In the same year, he made his first trip to Japan, where he was invited to the conference ``Manifolds'' in Tokyo. Thom was also there, as a guest of honor. Norbert remembers that Thom was very unhappy because he was placed in a special hotel,  a very noble one, but isolated and far from all the other participants. During that conference, on one afternoon, Ichiro Tamura presented his students to Norbert and allowed them to ask questions to him, one by one. Among these students was Mutsuo Oka who soon after (in the summer of the same year) went to France to work with Norbert as a PhD supervisor.  Oka obtained his doctorate at Orsay in 1975. He was the first PhD student of Norbert and a good friend ever since.

After having defended his thesis at Orsay, Norbert was hired Senior Researcher at CNRS. He had an office at Orsay and another one at IH\'ES (Bures). But he was uncomfortable with such a position. At CNRS, one is supposed to be a full-time researcher,  with no teaching duties. Norbert always liked the combination of research and teaching.  He stayed only one year at CNRS.
The year after (in 1974), he was appointed professor at Jussieu. In the same year, he  was an invited lecturer at the Vancouver ICM. He gave there a talk on the monodromy group of the unfolding of isolated singularities of planar curves. The topic had been the subject of an intensive activity since the beginning of the 1960s, with works of Hirzebruch, Arnold, Thom, Brieskorn, Tjurina, Milnor and others.
The participants of the ICM were staying on the campus of the university located next to the sea overlooking a very beautiful beach. One can imagine Norbert, after the lectures, spending the whole evening and part of the night talking mathematics, staring at the ocean and watching the rolling waves. He remembers that one evening, on that very beach, Douady organized a m\'echoui---a whole lamb grilled on a spit---for a large group of mathematicians.

Norbert taught at Jussieu for 3 years. 
In 1977, he exchanged his position  at Jussieu with that of Rosenberg, who was teaching at Orsay. Rosenberg used to live in Paris and he was happy with this arrangement. In this way he was not obliged anymore to go to Orsay several days a week. Symetrically, Norbert was happier in Orsay. It was the countryside. Furthermore, he preferred to be in Orsay because of the asbestos problem in Jussieu. Starting in 1974, people were talking about this problem there. Rosenberg said that this asbestos problem did not bother him.

   The famous Orsay seminar on Thurston's works on surfaces took place during the academic year 1976-77.  Norbert was not present; he spent that year in Lausanne.  But he used to return to Paris from time to time, and Laudenbach and others told him regularly what was going on at the seminar.

 Another occurrence at Orsay, involving Thurston, took place at that time.  People learned about the latter's work on the topology of 3-manifolds that uses hyperbolic geometry, from his Princeton lecture notes which reached Orsay, and whose  importance was immediately realized. The introduction of hyperbolic geometry in the study of the topology of 3-manifolds was completely new. Siebenmann asked Norbert to give a course on hyperbolic geometry. Norbert responded that he knew very little about this. Siebenmann told him then: ``This should not be a problem, you will learn it by giving this course." It was a happy coincidence that Norbert was planning to visit the Mittag-Leffler Institute, just after Siebenmann's request. There, in one of the attics that were next to the library, he 
found a set of old papers, notes and documents on hyperbolic geometry.  He went through them and he took some notes. He came back to France with enough material to build a course.
 
At Orsay, Douady and Po\'enaru knew a little bit of hyperbolic geometry. Deligne also had already come across this subject.  He was studying Hodge structures, which can be viewed as symmetric spaces, and he used to make computations using models of hyperbolic geometry.  
For most of the other mathematicians at Orsay, hyperbolic geometry was unknown at that time. It was a dormant subject and it needed to be resurrected.  
In fact, before Thurston arrived at Orsay, even the subject of Riemann surfaces in the original sense of Riemann, Klein and Poincar\'e---the sense which was restored by Thurston---, was very poorly known there. Algebraic geometry was at its peak. A Riemann surface was considered as a 1-dimensional scheme. Thurston brought hyperbolic geometry and Riemann surfaces at the center of the discussions at the topology seminar at Orsay.

Norbert gave his course, and as Siebenmann had predicted, he learned hyperbolic geometry while he was lecturing on this subject.
In his teaching, and later on, Norbert was doing hyperbolic geometry as it would have been done in school, or rather, as Euclid's geometry was taught in school: starting with the fact that by any two points passes a single line, and then the other axioms and postulates, then, explaining the geometry of triangles, showing that the three altitudes of any triangle intersect in one point, etc. Norbert also established the trigonometric formulae, and then he passed to dimension 3. His point of view was that of Lobachevsky, that is, model-free,

About twenty-five years after Norbert gave this course, I went with him through his notes. We expanded them together, we included spherical geometry and transitional geometries into the discussion and we published them. This is how I learned much of what I know in non-Euclidean geometry.

Norbert was professor at Orsay for the period 1977-1982. He taught differential calculus to second- and third-year students. The differential calculus course at French universities  has very little to do with the calculus courses given at American universities. There were a few textbooks available on the subject,  by Dieudonn\'e, Dixmier and some others, but it was (and is still) unusual for a professor in France to assign a book  or follow a book written by someone else. The students usually take notes. In 1977-78, Norbert's course was centered on the qualitative behavior of integral curves of differential equations. At least this is what the author of these lines, who was a third-year student there, remembers. The third-year students were divided into two groups. Douady taught the other group. Norbert recalls that some students and teaching assistants were complaining because he was doing an unusual program.

In 1979, Norbert gave a talk at the Bourbaki seminar. The subject was the first part of Hilbert's 16th problem, a question which concerns the number of connected components of a smooth algebraic plane curve in terms of the degree. He explained results of Gutkov, Arnol'd, Rokhlin, Guillou--Marin and others. I was a fourth-year student (in the ``ma\^\i trise" year) and I attended that seminar, just by curiosity. I  enjoyed the pictures, and the way Norbert was drawing them on the blackboard. Later on, I remembered these pictures when I saw for the first time the Chladni figures that appear in music theory, namely, in the theory vibration of plates.

During the period he spent at Orsay and Jussieu, Norbert had the opportunity to visit again the United States.   
In those times, it was possible for a university professor in France to do the total amount of his yearly teaching hours in one semester, and then, during the other semester, he could travel. Nobody complained about this. In 1977, David Eisenbud invited him to Brandeis. He gave there a course on singularities. 
Needless to say, traveling to the United States was  the occasion for him to meet mathematicians. Among those he met, he mentions Bott and Thurston.
 
 Back to Sullivan.

Dennis arrived at Orsay in 1973-74 as a visiting professor.
During that year, he conducted an explosive seminar at Bures, which gathered all the dynamicists from the Paris area. The sessions were always followed by long discussions, which took place first at  l'Ormaille and continued at Dennis's apartment, Boulevard Jourdan. This was an apartment given to Dennis by the university.  In fact, the whole building belonged to the City of Paris and it was managed by the university, and several mathematicians lived there for some time. Henri  Cartan used to live in that building, and he was coordinating  the distribution of the various apartments. Hadamard also had lived in that building. Pierre Lelong had occupied Dennis's apartment before him, and he had left it because he found it too big for him. Norbert remembers that officially, Lelong was still the tenant. Dennis liked to invite mathematicians at his home. During his visits to him, Norbert met several people. He remembers having seen there  Andr\'e Weil's daughter. 

After that visit to Orsay, Dennis was proposed a permanent position at IH\'ES; it was the position that was left vacant by Grothendieck.

Soon after, Norbert got acquainted with Swiss universities. He visited Geneva, where there was a strong mathematical department, with Haefliger, whom he already knew, Michel Kervaire, who was French but who had studied in Z\"urich, and who had decided to work in Switzerland  after a few years of teaching at New York University, and there were also other very good topologists and geometers. It was Kervaire who told Norbert about this possibility of visiting position.

In Geneva, there was a mathematics seminar where everyone could ask many questions, as many as he wanted. The participants' aim there was to understand.  ``Time was infinite there", Norbert says. After the seminar, the discussions  continued at the caf\'e. In contrast, asking questions during a seminar was rather unusual in France.

 An opportunity arose for Norbert to spend a year in Lausanne, which is less than one hour driving or by train from Geneva. He  took a leave from Paris and spent the academic year 1976-77 there.

In 1982, after having been professor at Orsay for five years (1977-78 to 1981-82), Norbert  left permanently France for Switzerland.  About this move to Switzerland, I remember the following facts I learned from him.

 During one of his visits to the Geneva mathematics section, Norbert saw a small handwritten announcement on the bulletin board, for a professorship  in Basel. He asked people there  for more information and for advice. The only response he heard was that the teaching in Basel was in German. In Geneva, very few people spoke German. An exception was Kervaire, who was fluent in several languages, including German and modern Greek. Haefliger told Norbert that for the Swiss living in French speaking regions, the German language reminds them of their military service. Norbert said he had no problem with teaching in German.  
 
He applied and obtained the position. He was appointed professor in Basel, in 1982. 
 
 I do not remember the exact year where Norbert started coming regularly from Basel to Strasbourg, but it was in the 1990s. He used to come by bike. Sometimes, when there was a conference, he would come by bike together with his students, but it would take longer, because the students were not used to such a long bike trip. Everyone who knows Norbert knows his passion for cycling. He can immediately recognize a good bike. I noticed that when he sees one on the street he stops to watch.  He told me that he bought his first good bike with his first professor salary, and it costed him the entire salary. People close to him were surprised by the price he spent on the bike. After a few years, this bike was damaged in an accident and he bought an even better one, which he still uses today, after more than 45 years.  
 
In May 2013, I was at a master-class on the island of Samos, and Norbert was among the teachers.  In order to be able to explore the island by bike, he took with him his mountain bike, on the plane. He needed to pump air after the flight and therefore entered a bike shop in Karlovassi. He still remembers the 
 owner's name, Georgios Katsouris. The latter gave him a pump and asked him to sign a paper written in Greek. He then revealed to him that with this signature he had registered for a bike race a few days later in a village called Kyriaki, in the mountains of Samos. Norbert indeed took part in this race and ended up first in the category of over 40s. After the race the entire village got together for a feast on the market place with a lot of food and drink. Norbert spent the evening there, talking to bikers and villagers, asking questions about the life in the island, enquiring about meanings of names and of words, while the other mathematicians were having their usual dinner at the conference hotel. This is how Norbert made friends in Samos. 
At the end of his stay, he had to pass by the store of Katsouris to return the pump, and I went with him. Katsouris used to sell and rent bicycles which mere made by his hands. Norbert was in awe: he told me that these hand-made bikes were really great.  He discussed with Katsouris for hours. I already knew (from Norbert) that there is a worldwide community of cyclists, and that whenever they meet they have millions of things to talk about.
Our common friend Charalampos Charitos tells me that some Greek mathematicians, when they talk about Norbert, talk about the ``cyclist". One day, we were at a conference in Kunming, and there is a big lake there. Charitos, who was also there, told Norbert that he would like to go around the lake by bicycle. Norbert took this very much to heart and spent some time finding him a bike. Eventually, Charitos realized that the distance to go around the lake was too long and he gave up. 

Norbert rides his bike mostly in the countryside. He crossed France by bicycle. I can understand that this passion for cycling is part of his need for freedom and open spaces. There is also, for a man who, like him, has a peasant's soul, the simple joy of contact with nature, a pleasure for the eyes and for the heart. To see the farms and all these fields around him reminds him of his childhood.  Once, I was at his place in Witterswil, and he told me: let's go to buy some bread. I realized he wanted us to go by bike and he gave me one of his bikes.  I thought we were just going to the next village, but in fact his plan was to go quite far. We took small country roads, sometimes crossing meadows and going around ponds.  He was riding his bike ahead of me, and the distance between us was getting bigger and bigger, so very soon I lost sight of him. At a crossroads I didn't know which way to go. As he doesn't have a cell phone, and since I didn't know the way back to his home, I waited for him at the crossroads. The time seemed long, I was anxious, but eventually he came back. When people go to his place, I tell them: beware, do not go for biking with him unless you are used to that.

 I asked Norbert why he left France for Switzerland. His first answer was ``teaching". He has always enjoyed teaching at all levels.    At Orsay,  he did not like the way teaching was distributed. 
  The first and second year students were not sufficiently taken care of. Students who finished the ``classes pr\'eparatoires aux grandes \'ecoles" entered directly the ``licence" year (third year) and the level  then rose sharply. The first year students of the \'Ecole Normale Sup\'erieure arrived also in that licence year.
These students received special instruction at their school, in addition to the courses they followed at Orsay, where the program seemed easy for them; they were only required to obtain their ``licence" and then their ``ma\^\i trise" there, since the \'Ecole Normale does not deliver such diplomas. They were being prepared for a PhD program.  
The result is that the third and fourth year students at Orsay received much more care than those of the first and second years. In general, most of the attention was devoted to those who had finished a ``grande \'ecole" (generally, \'Ecole Normale Sup\'erieure or \'Ecole polytechnique) and who came to Orsay as assistants or with support from CNRS to work on a PhD or to take part in seminars. All this was at the expense of the first and second year students who were being left behind. 
This system seemed absurd to Norbert.
 
Another reason for which Norbert preferred Switzerland is that in France, teachers were bound to a ``program", whereas in Switzerland they had more freedom.  
At Orsay, at that time, the teaching for first and second year students took place in parallel classes, each of them addressed to more than 200 students (sometimes much more) packed in big lecture halls, whereas Norbert was dreaming of teaching  courses in the spirit of those he had with Freudenthal and his other teachers in Utrecht. It was only later, in Basel, that he was able to have comparable conditions, and with no predefined program: he always made his teaching adapted to the students he had in front of him.

I  also asked Norbert how did it come about that he became president of the Swiss Mathematical Society. Being president of whatever organism does not fit with his personality. He responded: ``It  was a big  mistake". He did not really explain what he meant by that, but he recounted the circumstances under which he became president, and indeed it was a chain of misunderstandings. It started a day when he decided to go to the annual meeting of the Society, where Thom was announced as a speaker.
The annual meeting was organized that year in the new canton of Jura, at two different locations, namely Porrentruy and Del\'emont. Norbert did not read correctly the announcement, and on that day he went to Del\'emont, whereas Thom was speaking in Porrentruy.  Incidentally, Thom's hometown, Montb\'eliard (in France), is about 20 km far from  Porrentruy. There were very few people at Del\'emont where Norbert arrived; on that day, only the administrative meeting took place there. The mathematicians in the executive committee of the SMS chatted with him; they told him that a position of treasurer of the Society was vacant, and they convinced him to run for it. This is how he was elected treasurer. 
But soon after, he learned that normally the treasurer is elected secretary after two years, and then president for another two years. So eventually he was elected secretary then president. This is how Norbert became president of  the Swiss Mathematical Society, for the years 1988 and 1989. Norbert is not talented for administrative duties. He told me that one of the former presidents, his  friend and colleague in Basel Heinz Huber, when he learned about his election,   said: ``Let us see if the Society will survive that president."
 But he managed, even though he did not know much about Swiss organisation: the confusion between Porrentruy and Del\'emont is just one example of the things with which he was unfamiliar.


%
%
%
%

Recently, Norbert wrote a book on geometry, very different from all the existing ones.  The book is titled \emph{Topological, differential and conformal geometry of surfaces}. It is based on courses he gave in Basel, Vienna, Samos, Kunming, Bangalore and Varanasi. The book contains classical material that every student must have seen at least once: Complex structures from the point of view of $J$-fields, integrability on surfaces and the obstructions to integrability in higher dimensions by the Nijenhuis tensor, the Poincar\'e lemma, Brouwer's theorem, Morse functions, Stokes theorem, de Rham,  \v{C}ech and Dolbeault cohomologies, the uniformization theorem for compact Riemann surfaces, the construction of the Riemann surface associated with a meromorphic function, Riemann--Roch, the construction of Teichm\"uller's universal curve, the Weil--Petersson K\"ahler structure on Teichm\"uller space,  the embedding of Riemann surfaces in projective spaces and Chow's theorem, etc. But there are also a lot of special topics dear to Norbert:   models of the field $\mathbb{C}$ of complex numbers and its group of automorphisms, constructions with finite fields including several ones using the field with one element, a model for the algebraic closure of $\mathbb{C}$, 
the hyperbolic plane seen as the set of surjective homomorphisms of the ring $\mathbb{R}[X]$ of one-variable polynomials in a field modulo the automorphisms of the field, other models of the same plane: as the set of complex linear structures on the real plane, as the space of cross-ratio preserving involutions of the real projective line, etc. The book also contains several new models of hyperbolic 3-space.
  
  All along the book, when he mentions for the first time the name of a mathematician, Norbert recalls that this mathematician is Dutch, French, Iranian, British, Swiss, Russian, Polish, etc., exemplifying the fact that mathematics is the result of a collaboration without borders.

 We have come to the end of this short biography. I would like to take this opportunity to include a few more personal thoughts.

My relation with Norbert transformed during the years from that of a  student with his teacher, into friendship and collaboration. 
We published a certain number of joint papers, but I think that for me, from the purely mathematical point of view, the most rewarding collaboration with him was the fact that  at several occasions we went through writings of great authors from the past, including Euler, Riemann,  Klein, Teichm\"uller, Poincar\'e, Grothendieck and a few others. A priori it was for no special purpose, just for the satisfaction of reading the masters, but this turned out to have an impact on both of us. We understood important ideas there, and we used them later in our works.  
 
Norbert has a social and humanistic view on life, which dates back to his deeply agricultural origins.  He has that profound sense of justice, which also goes back to his peasant origins, that of people who know the honesty of the land that gives back exactly what it owes to the work of the farmer. As a mathematician, he carries a sense of intellectual honesty to a very high degree. He is a hard worker, but he also knows how to have moments of rest free of all worries. After a few hours of mathematical discussions, he tells me ``Now we're going to have a drink". I learned from him the sentence ``\`A la belle vie !" which he pronounces as an expression of good wishes.
He is curious about everything. He often asks questions about etymology and the pronunciation and spelling of words in foreign languages (Greek, Latin, Arabic, Chinese, Japanese, etc.), looking for relations between them. He is always curious about customs and ways of life in foreign countries. It is always interesting to hear him talking about scenes or small details which nobody else would have noticed, after a stay he has made at a mathematical institute in Pakistan, Algeria, etc.  Without putting it in appearance, he knows a lot about history and philosophy. He has a very unconventional and anti-establishment attitude about life and relations between people, but always with an extremely civic spirit. He has a refractory and insubordinate character. In his teaching at university, he did not want to hear about programs.  He hates administration and power. Without being an anarchist, he is, in every sense of the word,  the opposite of a bourgeois.   
Spending time with him gives you a real sense of freedom.

Their house, with Annette, is permanently open. His former (graduate and undergraduate) students often visit them.
 Everyone who knows that house can go there any time.  Common friends, when they come to Strasbourg, have taken the habit of visiting them in Basel; the distance between the two cities is only 120 km. I cannot mention them all, but recently there were Sumio Yamada, Ken'ichi Ohshika, Bob Penner, Alexey Sossinsky, Krishnendu Gongopadhyay, and there are many others. To his visitors who come for the first time to Basel, Norbert is used to show the tombstone of Jacob Bernoulli in one of the cloisters of Basel's cathedral, decorated with the logarithmic spiral design and the motto \emph{Eadem mutata resurgo} (``Changed and yet the same, I rise again"). He also shows them the house where Euler lived as a child and the church where his father officiated as pastor, before taking them to his house,  in the middle of a garden, part of which he cultivates himself, and the other part cultivated by Annette.
  Friends have become collaborators and vice versa.
 
In his book on geometry I mentioned above, after a section in which he explains the classification of surfaces using Morse functions, he writes, in a small section entitled \emph{Thoughts}: ``In societies one tries to measure trends, for instance `optimistic', `pessimistic', `undecided' by $o-p$, and forgets $u$. Situations which invite the response `yes', `no' or `doubt' tend to be measured by $y-n$, forgetting the $d$, which is a counting that systematically forgets those who want to think or think more."
He then adds: ``Unfortunately, many countries on today's Earth consider individuals as critical points. This is true especially in the context of power and truth. In a society with only one source of truth, say by a religion or by a dominant, or only one, political party, the count is even more restricted. An individual with an opinion that does not question the official truth is counted, those that do question, but keep silent, are not counted, and others are stored in camps and prisons, if not eliminated."

I would like to end with a few words that Norbert wrote to me in a recent mail (November 1st, 2022): ``I see mathematics as an activity directly related to life, the life of animals, and therefore also of humans," a sentence which immediately reminded me of Manin's sentence, from his ICM talk \emph{Mathematics as a metaphor} (Kyoto, 1990): ``Mathematics is a novel about Nature and Human kind", and also, of  Galileo Galilei's often quoted related sentence:  ``[Natural] Philosophy is written in this very great book that continually stands open before our eyes (I say the universe), but it cannot be understood unless one first learns to understand the language, and the characters, and what is written there. It is written in a mathematical language, and the characters are triangles, circles, and other geometrical figures, without which it is impossible to humanly understand a word of it; without them, it is a vain wandering through a dark labyrinth." (The Assayer, 1623).

Norbert continues:
``The migration of birds, like the hunting methods of raptors, require an enormous geometrical knowledge.
Look at a bird landing directly on a tiny shaking branch: no artificial intelligence can yet achieve this. Similarly, the cow that
starts licking the newborn calf's skin, immediately after birth, in one direction and then in the other, so that the hairs form an air cushion which provides thermal protection,
is very competent in differential geometry." He then adds: ``But
the human species can communicate with future generations by writing, while other species do it rather by a selection that is inscribed in the genes, and this process is much longer. They cannot give a knowledge to the
generation that follows them immediately. This gives humans a slight advantage. Through writing, we can pass on direct knowledge to our offspring."   

In a technical civilization where we sometimes see people engaged in senseless activities, a friendship with someone like Norbert helps to regain a balanced view on life.

       %
%
%
%
%
%
%
%


\end{document}